\documentclass[11pt,a4paper,leqno]{article}
\usepackage{amsmath,amsthm,amsfonts,amssymb}
\newtheorem{defin}{Definition}[section]

\newtheorem{propos}[defin]{Proposition}
\newtheorem{theor}[defin]{Theorem}

\newenvironment{definition*}{\begin{defin} \rm}{\end{defin}}
\newcommand{\BOX}{\rule{2mm}{2mm}}
\newcommand{\Zint}{{\mathbb {Z}}}
\newcommand{\Qq}{\mathbb Q(q)}

\newcommand{\eN}{{}e}
\newcommand{\fN}{{}f}
\newcommand{\tN}{{}t}
\newcommand{\eL}{{}e}
\newcommand{\fL}{{}f}
\newcommand{\tL}{{}t}
\newcommand{\vN}{v^{(n)}}
\newcommand{\vL}{v^{(l)}}
\newcommand{\UN}{U_q^{\prime}(\widehat{{\mathfrak {sl}}}_n)}
\newcommand{\UL}{U_q^{\prime}(\widehat{{\mathfrak {sl}}}_l)}
\newcommand{\UNd}{U_q(\widehat{{\mathfrak {sl}}}(n))}
\newcommand{\UNdd}{U_q(\widehat{{\mathfrak {sl}}}_n)}
\newcommand{\VNL}{V_{n,l}}
\newcommand{\DeN}{{}^n\hspace{-0.8mm}\Delta}
\newcommand{\DeL}{{}^l\hspace{-0.8mm}\Delta}
\newcommand{\ba}{\boldsymbol {a}}
\newcommand{\bb}{\boldsymbol {b}}
\newcommand{\bh}{\boldsymbol {h}}
\newcommand{\bk}{\boldsymbol {k}}
\newcommand{\bll}{\boldsymbol {l}}
\newcommand{\akW}{{}^{\ba(\bk)}\hspace{-0.5mm}\widehat{W}}
\newcommand{\aW}{{}^{\ba}\hspace{-0.5mm}\widehat{W}}
\newcommand{\bs}{\boldsymbol {s}}
\newcommand{\bl}{\boldsymbol {\lambda}}
\newcommand{\bm}{\boldsymbol {\mu}}
\newcommand{\bo}{\boldsymbol {0}}
\newcommand{\FN}{ {}^n\hspace{-0.5mm}{\mathcal{F}}}
\newcommand{\FL}{ {}^l\hspace{-0.5mm}{\mathcal{F}}}
\renewcommand{\P}{\mathcal {P}}
\newcommand{\A}{\mathcal {A}}
\newcommand{\B}{\mathcal {B}}
\newcommand{\ANL}{A^{(n,l)}}
\newcommand{\ALN}{A^{(l,n)}}
\newcommand{\DNLp}{{}^+\hspace{-0.5mm}D^{(n,l)}}

\newcommand{\DNLm}{{}^-\hspace{-0.5mm}D^{(n,l)}}

\newcommand{\DNLpm}{{}^{\pm}\hspace{-0.5mm}D^{(n,l)}}
\newcommand{\DLNpm}{{}^{\pm}\hspace{-0.5mm}D^{(l,n)}}
\newcommand{\DNLmp}{{}^{\mp}\hspace{-0.5mm}D^{(n,l)}}

\newcommand{\lL}{\bl_l}
\newcommand{\lN}{\bl_n}
\newcommand{\sL}{\bs_l}
\newcommand{\sN}{\bs_n}
\newcommand{\mL}{\bm_l}
\newcommand{\mN}{\bm_n}
\newcommand{\msL}{ \bm_l, \bs_l }
\newcommand{\msN}{ \bm_n, \bs_n }
\newcommand{\lsL}{ \bl_l, \bs_l }
\newcommand{\lsN}{ \bl_n, \bs_n }
\newcommand{\ph}{\varphi}

\newcommand{\fra}[1]{
\begin{picture}(8,8)(2,1.5)
\multiput(0,0)(0,8){2}{\line(1,0){8}}
\multiput(0,0)(8,0){2}{\line(0,1){8}}
\put(0,0){\makebox(8,8){{\tiny $ #1 $}}}
\end{picture}}
\newcommand{\ov}{\overline}
\newcommand{\barr}{\,\ov{{\makebox(6,6){}}}\,}

\begin{document}
\thispagestyle{empty}
\begin{center}
%
%
{\LARGE 
Canonical bases of higher-level $q$-deformed \\[2mm] Fock spaces}
\\[8mm]
%
%
{\Large  Denis Uglov} \\[8mm]
%
%
\parbox{5.42in}
{We define canonical bases of the higher-level $q$-deformed Fock space modules
of the affine Lie algebra $\widehat{{\mathfrak {sl}}}_n$ generalizing the result of Leclerc and Thibon for the case of level 1. 
We  express the transition matrices between the canonical bases and the natural bases of the Fock spaces in terms of certain affine Kazhdan-Lusztig polynomials.}  
\end{center}
\medskip
\section{Introduction}
Leclerc and Thibon defined, in \cite{LT1}, a canonical basis of the level-1 $q$-deformed Fock space module \cite{KMS,S} of $\UNdd,$ and conjectured that the entries of the transition matrix between this basis and the natural basis of the Fock space are $q$-analogues of the decomposition numbers of the $v$-Schur algebra. Recently this conjecture has been proved by Varagnolo and Vasserot \cite{VV} who gave, in particular, expressions for the transition matrix in terms of Kazhdan-Lusztig polynomials. These developments are summarized in the recent review \cite{LT2}. 

The aim of this note is to define canonical bases for the $q$-deformed Fock spaces with arbitrary positive integral levels. 
The higher-level $q$-deformed Fock spaces were introduced in \cite{JMMO} where they were used to compute the crystal graph of an arbitrary irreducible integrable module of $\UNdd.$ 
More recently these spaces appeared in the work \cite{FLOTW} whose conventions we largely follow here. 

A $q$-deformed Fock space of level $l$, denoted $\FN_{\sL},$ is an integrable $\UN$-module parameterized by a sequence $\sL=(s_1,\dots,s_l)$ of integers, it has the natural basis $\{\ph(\lsL)\:|\:\lL\in \Pi^l\}$ where $\Pi^l$ is the set of $l$-multipartitions. 
The submodule $M_{\sL}$ of $\FN_{\sL},$ generated by the {\em highest weight vector} $\ph(\bo_l,\sL),$ where $\bo_l$ is the $l$-tuple of zero partitions, is isomorphic to the irreducible $\UN$-module with highest weight $ \Lambda_{s_1\bmod n} + \cdots + \Lambda_{s_l\bmod n}.$  

Our definition the canonical basis, denoted $\{ G^+(\lsL)\: | \: \lL \in \Pi^l\},$ of $\FN_{\sL}$ is quite similar to that of $\cite{LT1}.$ We use the semi-infinite $q$-wedge construction of the higher-level $q$-deformed Fock spaces developed in \cite{TU}. 
This construction  allows to define a natural semi-linear involution of the Fock space leaving invariant the highest weight vector and commuting with the lowering generators of $\UN.$ 
The unitriangularity of the involution, understood  in the sense described in Section 3.1, makes it possible to define the canonical basis in an elementary way. 
The canonical basis is compatible with the inclusion   of $M_{\sL}$ into $\FN_{\sL}:$ the elements $G^+(\lsL)$ indexed by multipartitions $\lL$ from the set $\Pi^l_{\sL}$ described in Section 3.2, coincide with the elements of the global lower crystal basis \cite{K} of $M_{\sL}.$  

An algorithm for computing the global lower crystal basis of an arbitrary irreducible integrable module of $\UNdd$ has been given in \cite{Ma}. 
Our construction gives a straightforward recursive procedure that allows to compute the entries of the transition matrix, denoted $\DNLp_{\lL,\mL}(\sL;q),$ between the canonical basis  $\{G^+(\lsL)\}$  and the natural basis $\{\ph(\lsL)\},$ thereby giving a different approach to the problem of $\cite{Ma}.$       

By the deep result of Ariki \cite{A} the value at $q=1$ of the transition matrix between the lower global crystal base of $M_{\sL}$ and the natural basis of $
\FN_{\sL}$ gives the decomposition matrix of a certain Ariki-Koike algebra. It would be interesting to find a similar interpretation for the canonical basis of the entire Fock space $\FN_{\sL}.$

\vspace{2mm}
\noindent{\bf {\Large Acknowledgments }}\\
I am grateful to T. Baker, M. Kashiwara, T. Miwa, K. Takemura and J.-Y. Thibon for discussions. Special thanks are due to Bernard Leclerc for many patient explanations concerning the subject of works \cite{FLOTW, LT1, LT2}.

\section{Higher-level $q$-deformed Fock spaces}
\subsection{ The $q$-wedge products}
Let $V_n$ be an $n$-dimensional vector space over $\Qq$ with the basis $\vN_1,\dots,\vN_n.$ Throughout this note we shall fix integers $n > 1$ and $l > 1,$ and put $V_{n,l} = (V_n \otimes  V_l)[z,z^{- 1}].$ Let $\UN$ be the quantum group corresponding to the Lie algebra $\widehat{{\mathfrak {sl}}}_n^{\prime}.$ This is an associative algebra over $\Qq$ with generators $\eN_i,\fN_i,\tN_i,\tN_i^{-1},$ $i=0,\dots,n-1$ and the standard relations which may be found, for example, in \cite{KMS}. We define on $V_{n,l}$ the structure of a representation of $\UN$ by      
\begin{eqnarray}
\tN_i(\vN_a z^m \vL_b) & = & q^{\delta(i\equiv a\bmod n) - \delta(i+1\equiv a\bmod n)} \;\vN_a z^m \vL_b, \label{eq:t}\\
\eN_i(\vN_a z^m \vL_b) & = & \delta(i+1\equiv a\bmod n) \;\vN_{a-1} z^{m + \delta(i=0)} \vL_b, \label{eq:e}\\
\fN_i(\vN_a z^m \vL_b) & = & \delta(i\equiv a\bmod n) \;\vN_{a+1} z^{m - \delta(i=0)} \vL_b, \label{eq:f}
\end{eqnarray}
where it is understood that $\vN_0 = \vN_n,$ $ \vN_{n+1} = \vN_1$ and for a statement $P$ we put  $\delta(P) =1$ if $P$ is true,  $\delta(P) =0$ if $P$ is false. 
We also define on $V_{n,l}$ the structure of a representation of $\UL$ by exchanging everywhere in the preceding formulas $n$ with $l$ and $a$ with $b.$ Clearly  the actions of $\UN$ and $\UL$ on $V_{n,l}$ obtained in this way commute one with another.

In the sequel it will be sometimes convenient to label elements of the basis $\{ \vN_a z^m \vL_b \: |\: a\in \{1,\dots,n\}, b\in \{1,\dots,l\}, m\in \Zint\}$ of $\VNL$ by a single integer : we put $k=a +n(b-1) - nl m $ and $u_k = \vN_a z^m \vL_b$ so that $\{ u_k \:|\: k\in \Zint \}$ is a basis of $\VNL.$     

The {\em $q$-wedge square } of  $\VNL,$ denoted $\bigwedge ^2 \VNL,$  is  a $q$-deformation of the exterior square of $\VNL.$ It is a  $\Qq$-vector space generated  by monomials  $u_{k_1}\wedge u_{k_2}$  $(k_i \in \Zint)$ which obey the defining relations described as follows. First of all, the monomials $u_{k_1}\wedge u_{k_2}$ with $k_1 > k_2$ form a basis of $\bigwedge ^2 \VNL,$ such monomials are called {\em ordered}. 
Any  monomial $u_{k_1}\wedge u_{k_2}$ with $k_1 \leq  k_2$ is expressed as a linear combination of these by means of the ordering rules (cf. \cite[Lemma 3.2]{TU} where $q$ is $q^{-1}$ in the present notations): put $k_i = a_i + n(b_i -1) - nlm_i$ where $a_i\in \{1,\dots,n\},$ $b_i \in \{1,\dots,l\},$ $m_i \in \Zint.$ Let $ \alpha = (a_2-a_1)\pmod {nl},$ $ \beta = (n(b_2-b_1))\pmod {nl},$ then for $k_1\leq k_2$ 
\begin{align}
&\text{if} \; \alpha=0,\beta=0:  \quad \quad u_{k_1}\wedge u_{k_2} = - u_{k_2}\wedge u_{k_1}; \tag{Ri} \label{eq:no1}\\
&\text{if} \; \alpha>0,\beta=0: \quad \quad u_{k_1}\wedge u_{k_2} =  - q^{-1} u_{k_2}\wedge u_{k_1} + \tag{Rii}  \\
  & + (q^{-2}-1)\{\sum_{m\geq 0} q^{-2m} u_{k_2 - \alpha -nlm} \wedge u_{k_1 + \alpha + nlm}\; - \; \sum_{m\geq 1} q^{-2m+1} u_{k_2 - nlm} \wedge u_{k_1 + nlm} \}; \nonumber  \\
&\text{if} \; \alpha=0,\beta>0: \quad \quad u_{k_1}\wedge u_{k_2} =  - q u_{k_2}\wedge u_{k_1} + \tag{Riii}  \\
  & + (q^2-1)\{\sum_{m\geq 0} q^{2m} u_{k_2 - \beta -nlm} \wedge u_{k_1 + \beta + nlm}\; - \;\sum_{m\geq 1} q^{2m-1} u_{k_2 - nlm} \wedge u_{k_1 + nlm}\};\nonumber   
\end{align}\begin{align}
&\text{if} \; \alpha>0,\beta>0: \quad \quad u_{k_1}\wedge u_{k_2} =  -  u_{k_2}\wedge u_{k_1} + \tag{Riv} \label{eq:no4} \\
&+(q-q^{-1})\sum_{m\geq 0}\frac{(q^{2m+1}+q^{-2m-1})}{(q+q^{-1})} \;
u_{k_2 - \beta -nlm} \wedge u_{k_1 + \beta + nlm}  - \nonumber \\
& -(q-q^{-1})\sum_{m\geq 0}\frac{(q^{2m+1}+q^{-2m-1})}{(q+q^{-1})} \;
u_{k_2 - \alpha -nlm} \wedge u_{k_1 + \alpha + nlm} + \nonumber \\   
&+(q-q^{-1})\sum_{m\geq 1}\frac{(q^{2m}-q^{-2m})}{(q+q^{-1})} \; 
u_{k_2 +nl - \alpha - \beta -nlm} \wedge u_{k_1 -nl +\alpha + \beta + nlm} - \nonumber \\   
&-(q-q^{-1})\sum_{m\geq 1}\frac{(q^{2m}-q^{-2m})}{(q+q^{-1})}  \; 
u_{k_2  -nlm} \wedge u_{k_1 + nlm},  \nonumber 
\end{align}
where the summations continue as long as the monomials appearing under the sums  are ordered. 
Relations (\ref{eq:no1} -- \ref{eq:no4}) where $k_1,k_2$ range over all pairs such that $k_1 \leq k_2$ constitute the complete set of defining relations in $\bigwedge^2 \VNL.$  

For any integer $r\geq 2$ the $r$-fold $q$-wedge product $\bigwedge^r \VNL$ is the vector space generated by elements $u_{k_1}\wedge u_{k_2} \wedge \cdots \wedge u_{k_r}$ $(k_i \in \Zint)$ modulo the relations (\ref{eq:no1} -- \ref{eq:no4}) in every adjacent pair of the factors. 
The ordered monomials, i.e. $u_{k_1}\wedge u_{k_2} \wedge \cdots \wedge u_{k_r}$ with $k_1 > k_2 > \dots >k_r$ form a basis of  $\bigwedge^r \VNL.$  

Iterating $r-1$ times the coproducts $\DeN : \UN \rightarrow \UN^{\otimes 2}$ and $\DeL : \UL \rightarrow \UL^{\otimes 2}$ given by 
\begin{alignat}{3}
\DeN(\eN_i) &= \eN_i \otimes \tN_i^{-1} + 1\otimes  \eN_i,&\quad \DeN(\fN_i) &= \fN_i\otimes 1 + \tN_i \otimes \fN_i,&\quad \DeN(\tN_i) &= \tN_i\otimes \tN_i, \label{eq:copn}  \\
\DeL(\eL_i)& = \eL_i\otimes \tL_i + 1\otimes \eL_i,&\quad \DeL(\fL_i) &= \fL_i\otimes 1 + \tL_i^{-1} \otimes \fL_i,&\quad \DeL(\tL_i) &= \tL_i\otimes \tL_i, \label{eq:copl}
\end{alignat}
we define on $\bigwedge^r \VNL$  a structure of a $\UN \otimes \UL$-module. Note that the coproducts are to be chosen differently for $\UN$ and $\UL$ in order to ensure compatibility of the actions with the relations in  $\bigwedge^r \VNL.$  

For $s\in \Zint$ we define the {\em semi-infinite $q$-wedge product of charge $s$}, denoted $\bigwedge^{s+\frac{\infty}{2}}\VNL,$ as the inductive limit of  $\bigwedge^r \VNL$ where the maps $\bigwedge^r \VNL \rightarrow \bigwedge^{r+1} \VNL $ are given by $v \mapsto v\wedge u_{s-r}.$ 
For $v\in \bigwedge^r \VNL$ let  $v\wedge u_{s-r}\wedge u_{s-r-1}\wedge \cdots$ stand for the image of $v$ in  $\bigwedge^{s+\frac{\infty}{2}}\VNL,$ then  $\bigwedge^{s+\frac{\infty}{2}}\VNL$ can be thought of as the vector space generated by semi-infinite monomials $u_{k_1}\wedge u_{k_2}\wedge \cdots $ $(k_i \in \Zint)$ such that $k_i=s-i+1$ for $i\gg 1.$ 
Such a monomial is said to be ordered if $k_1 > k_2 > \cdots.$ As in the finite situation the ordered monomials form a basis of $\bigwedge^{s+\frac{\infty}{2}}\VNL.$ We let $|s)$ to stand for the distinguished vector $u_s\wedge u_{s-1}\wedge \cdots $ of $\bigwedge^{s+\frac{\infty}{2}}\VNL.$     

The semi-infinite $q$-wedge product may be endowed with the structure of $\UN\otimes \UL$-module by, roughly speaking, iterating the coproducts (\ref{eq:copn}) and (\ref{eq:copl}) infinite number of times. 
More precisely, the action of $\UN$ say, is defined as follows. For $u=u_{k_1}\wedge u_{k_2}\wedge \cdots $ $ \in $ $  \bigwedge^{s+\frac{\infty}{2}}\VNL$ set  
\begin{equation} \label{eq:f(u)}
\fN_i(u) = \fN_i(u_{k_1})\wedge u_{k_2}\wedge u_{k_3} \wedge \cdots + 
 \tN_i(u_{k_1})\wedge \fN_i(u_{k_2})\wedge u_{k_3} \wedge \cdots +  \: \dots \:. 
\end{equation}
Using the ordering rules (\ref{eq:no1} -- \ref{eq:no4}) and the condition $k_{i+1} = k_i -1 $ for $i\gg 1$ one verifies that the sum in the right-hand side of (\ref{eq:f(u)}) contains only a finite number of non-vanishing terms, hence $\fN_i$ is a well-defined operator on  $  \bigwedge^{s+\frac{\infty}{2}}\VNL.$ 
The action of the Cartan part of $\UN$ is defined as follows. 
Every $u\in \bigwedge^{s+\frac{\infty}{2}}\VNL$ is represented as $u = v\wedge |nlm)$ for a suitable (non-unique) choice of $m\in \Zint$ and $v\in \bigwedge^{r}\VNL.$ Put
\begin{equation} \label{eq:t(u)}
\tN_i(|nlm)) = q^{l\delta(i=0)} |nlm) \quad \text{and} \quad \tN_i(v\wedge |nlm)) = \tN_i(v)\wedge\tN_i(|nlm)).
\end{equation}
Using (\ref{eq:t}) one verifies that $\tN_i$ does not depend on the particular choice of $m$ and $v.$ 
To define the action of the raising generators $\eN_i,$ represent $u\in \bigwedge^{s+\frac{\infty}{2}}\VNL$ as $u = v\wedge |m)$ for a suitable (non-unique) choice of $m\in \Zint$ and $v\in \bigwedge^{r}\VNL,$ and put    
\begin{equation} \label{eq:e(u)}
\eN_i(u) = \eN_i(v)\wedge \tN_i^{-1}(|m)).
\end{equation}
The ordering rules (\ref{eq:no1} -- \ref{eq:no4}) are used to show that $\eN_i$ is a well-defined operator on $\bigwedge^{s+\frac{\infty}{2}}\VNL.$ 

The action of $\UL$ on $\bigwedge^{s+\frac{\infty}{2}}\VNL$ is defined similarly. One only has to exchange $n$ with $l$ everywhere in the formulas (\ref{eq:f(u)} -- \ref{eq:e(u)}) and taking  account of the different coproduct for $\UL,$ replace $t$ by $t^{-1}$ in $(\ref{eq:f(u)})$ and  $(\ref{eq:e(u)}).$ 
Below we describe the action of $\UN$ on $\bigwedge^{s+\frac{\infty}{2}}\VNL$ in a combinatorial way. For now note that (\ref{eq:t(u)}) implies that the level of the $\UN$-action is $l,$ and the level of the $\UL$-action is $n.$  

For non-zero $m\in \Zint$ we define an operator $B_m :  \bigwedge^{s+\frac{\infty}{2}}\VNL \rightarrow \bigwedge^{s+\frac{\infty}{2}}\VNL$ by setting for $u = u_{k_1}\wedge u_{k_2} \wedge \cdots $ $:$ 
\begin{equation}
B_m(u) = u_{k_1 -nlm}\wedge u_{k_2} \wedge \cdots  + u_{k_1}\wedge u_{k_2 -nlm} \wedge \cdots + \dots \:. 
\end{equation}
As in (\ref{eq:f(u)}) the ordering rules imply that the sum in the right-hand side contains a finite number of non-zero terms, hence $B_m$ is well defined. 
By construction the operators $B_m$ commute with the action of $\UN\otimes\UL.$ It was shown in \cite{TU} that they generate the Heisenberg algebra ${\mathcal H}$, i.e. satisfy the relations   
\begin{displaymath}
[B_k,B_m] = \delta(k+m=0) \: \gamma_k(q),
\end{displaymath}
where $\gamma_k(q)$ is a non-zero element of $\Zint[q,q^{-1}].$ 

The irreducible decomposition of $\bigwedge^{s+\frac{\infty}{2}}\VNL$ as $\UN\otimes\UL \otimes {\mathcal H}$ -- module is given in \cite[Theorem 4.10]{TU}.

\subsection{The $q$-deformed Fock spaces}
It will be convenient to introduce three different labelings for a semi-infinite ordered monomial 
$u_{k_1}\wedge u_{k_2} \wedge u_{k_3}\wedge \cdots $ $ \in $ $ \bigwedge^{s+\frac{\infty}{2}}\VNL$.  

{\bf 1)} Labeling by a pair $(\lambda , s)$ where $\lambda =(\lambda_1,\lambda_2,\dots )$ is a partition and $s\in \Zint,$ is obtained by defining  $s$ to be the charge of $\bigwedge^{s+\frac{\infty}{2}}\VNL,$ and by setting $\lambda_i = k_i - s + i -1.$   

{\bf 2)} Labeling by a pair $(\bl_l,\bs_l),$ where $\bl_l = (\lambda^{(1)},\dots,\lambda^{(l)})$ is an $l$-tuple of partitions and $\bs_l = (s_1,\dots,s_l) \in \Zint^l,$ is obtained as follows. Write $k_i = a_i + n(b_i -1) -nl m_i$ $(a_i \in \{1,\dots,n\}, b_i\in \{1,\dots,l\}, m_i \in \Zint).$ 
For each $b\in \{1,\dots,l\}$ let $(k_1^{(b)} > k_2^{(b)} > \dots )$ be the semi-infinite sequence obtained by arranging elements of the set $\{ a_i - n m_i\:|\: b_i = b\}$ in strictly decreasing order. 
Then there is a unique $s_b \in \Zint$ such that $k_i^{(b)} = s_b -i+1$ for $i\gg 1.$  Define the partition $\lambda^{(b)} = (\lambda^{(b)}_1,\lambda^{(b)}_2,\dots)$ by $\lambda^{(b)}_i = k_i^{(b)} - s_b + i -1,$ and set   $\bs_l = (s_1,\dots,s_l).$ Observe that  $\sum_{b=1}^l s_b  = s.$  

{\bf 3)} Labeling by a pair $(\bl_n,\bs_n),$ where $\bl_n = (\lambda^{(1)},\dots,\lambda^{(n)})$ is an $n$-tuple of partitions and $\bs_n = (s_1,\dots,s_n) \in \Zint^n,$ is obtained similarly. 
Keeping notations as in {\bf 2)}, for each $a\in \{1,\dots,n\}$ let $(k_1^{(a)} > k_2^{(a)} > \dots )$ be the semi-infinite sequence obtained by arranging elements of the set $\{ b_i - l m_i\:|\: a_i = a\}$ in strictly decreasing order. Then there is a unique $s_a \in \Zint$ such that $k_i^{(a)} = s_a -i+1$ for $i\gg 1.$  
Define the partition $\lambda^{(a)} = (\lambda^{(a)}_1,\lambda^{(a)}_2,\dots)$ by $\lambda^{(a)}_i = k_i^{(a)} - s_a + i -1,$ and set   $\bs_n = (s_1,\dots,s_n).$ Observe that  $\sum_{a=1}^n s_a  = s.$

Let $\Pi$ be the set of all partitions. Define maps
\begin{displaymath}
\iota_l : \Pi \times \Zint \rightarrow \Pi^l\times \Zint^l \qquad \text{and} \qquad \iota_n : \Pi \times \Zint \rightarrow \Pi^n\times \Zint^n 
\end{displaymath}
by setting in the notations of {\bf 1)} -- {\bf 3)} :  
\begin{displaymath}
\iota_l : (\lambda , s) \mapsto (\bl_l,\bs_l) \qquad \text{and} \qquad\iota_n : (\lambda , s) \mapsto (\bl_n,\bs_n).
\end{displaymath}
It is not difficult to see that $\iota_l$ and $\iota_n$ are bijections. Indeed, with $s$ fixed, $\bl_n$ is an $n$-quotient of $\lambda,$ and $\bs_n$ encodes the $n$-core of $\lambda,$ while  $\bl_l$ and  $\bs_l$ may be thought of as variants of an $l$-quotient and the $l$-core \cite[Ch.1, \S 1]{Mac}.   

For a semi-infinite strictly decreasing sequence $k_1,k_2,\dots$ such that $k_i=s-i+1$ $(i\gg 1)$ and the partition $\lambda=(\lambda_1,\lambda_2,\dots)$ defined by  $k_i = \lambda_i + s-i+1$ we put 
\begin{alignat*}{4}
&|\lambda, s\rangle & & = & &\quad u_{k_1}\wedge u_{k_2} \wedge \cdots \;,&  &{} \\
&| \bl_l, \bs_l \rangle & &= & &\quad |\lambda, s\rangle \qquad & &\text{if $(\bl_l, \bs_l) = \iota_l(\lambda,s),$} \\  
&| \bl_n, \bs_n \rangle & & = & &\quad |\lambda, s\rangle \qquad & &\text{if $(\bl_n, \bs_n) = \iota_n(\lambda,s).$}
\end{alignat*}

For $\bs_l  = (s_1,\dots,s_l) \in \Zint^l$ and $\bs_n = (s_1,\dots,s_n )\in \Zint^n$ we define  
\begin{equation*}
 \FN_{\bs_l} :=   \oplus_{\bl_l \in \Pi^l} \Qq | \bl_l, \bs_l \rangle, \quad \text{and}\quad   
\FL_{\bs_n} :=   \oplus_{\bl_n \in \Pi^n} \Qq | \bl_n, \bs_n \rangle
\end{equation*}
so that 
\begin{eqnarray}
\bigwedge^{s+\frac{\infty}{2}} \VNL & = & \bigoplus_{s_1+\cdots + s_l = s} \: \FN_{\bs_l},  \label{eq:weightdecL}\\ 
\bigwedge^{s+\frac{\infty}{2}} \VNL & = & \bigoplus_{s_1+\cdots + s_n = s} \: \FL_{\bs_n}. \label{eq:weightdecN} 
\end{eqnarray}
It follows from (\ref{eq:t(u)}) that $| \bl_l, \bs_l \rangle$ is a weight vector of $\UL$ of weight $(l+s_l -s_1)\Lambda_0 + (s_1-s_2) \Lambda_1 + \cdots + (s_{l-1} - s_l) \Lambda_{l-1},$ and that $| \bl_n, \bs_n \rangle$ is a weight vector of $\UN$ of weight $(n+s_n -s_1)\Lambda_0 + (s_1-s_2) \Lambda_1 + \cdots + (s_{n-1} - s_n) \Lambda_{n-1}.$ 
Therefore (\ref{eq:weightdecL}) represents the weight decomposition of $\bigwedge^{s+\frac{\infty}{2}}\VNL$ with respect to the action of $\UL,$ while (\ref{eq:weightdecN}) represents the weight decomposition of $\bigwedge^{s+\frac{\infty}{2}}\VNL$ with respect to the action of $\UN.$ 

Since the action of $\UN$ commutes with that of $\UL,$ each of the subspaces $\FN_{\bs_l}$ is stable with respect to $\UN,$ and  each of the subspaces $\FL_{\bs_n}$ is stable with respect to $\UL.$ 

The action of $\UN$ on $\FN_{\bs_l}$ may be presented combinatorially. To this end we  introduce a new basis $\{ \ph(\lambda,s)\:|\: \lambda \in \Pi\}$ of $\bigwedge^{s+\frac{\infty}{2}}\VNL.$ 
Elements of this basis differ from the vectors $|\lambda, s\rangle $ only by signs. 

For $(\lsL) \in \Pi^l\times \Zint^l$ we define $\Phi(\lsL)\in \Zint$ as follows. Let $\bo_l$ denote the $l$-tuple of zero partitions, and let   
\begin{align*}
&|\lsL \rangle  = u_{k_1}\wedge u_{k_2} \wedge \cdots \quad \text{ where $k_1 > k_2 > \cdots ,$}\\
&|\bo_l,\sL \rangle  = u_{k^{\circ}_1}\wedge u_{k^{\circ}_2} \wedge \cdots \quad \text{ where $k^{\circ}_1 > k^{\circ}_2 > \cdots .$}
\end{align*}
Writing, as in {\bf 2)} at the beginning of this subsection, $k_i=a_i+n(b_i-1) -nlm_i,$ $k^{\circ}_i=a^{\circ}_i+n(b^{\circ}_i-1) -nlm^{\circ}_i,$ set for $j=1,2,\dots$ 
\begin{displaymath}
\Phi_j(\lsL) = \#\{ i< j \: |\: b_i < b_j\} - \#\{ i< j \: |\: b^{\circ}_i < b^{\circ}_j\}. 
\end{displaymath}
Since  $|\lsL \rangle $ and $|\bo_l,\sL \rangle$ both belong to $\FN_{\bs_l},$ hence have the same weight with respect to $\UL,$ we have $\Phi_j(\lsL) =0$ for all large enough $j.$ We set $\Phi(\lsL) = \sum_{j=1}^{\infty} \Phi_j(\lsL),$ and define  
\begin{alignat}{4}
&\ph(\lsL) & &:= & &\quad(-1)^{\Phi(\lsL)} \, |\lsL \rangle , \label{eq:ph1} & &\\
&\ph(\lambda,s) & & := & &\quad \ph(\lsL)\qquad & &\text{if $ (\lsL) =\iota_l(\lambda,s),$} \label{eq:ph2}\\
&\ph(\lsN) & & := & &\quad \ph(\lambda,s)\qquad & & \text{if $ (\lsN) =\iota_n(\lambda,s).$} \label{eq:ph3}
\end{alignat} 
Thus $\{ \ph(\lambda,s) \: | \: \lambda \in \Pi\}$ is a basis of $\bigwedge^{s+\frac{\infty}{2}}\VNL,$  $\{\ph(\lsL)\:|\: \lL \in \Pi^l\}$ is a basis  of $\FN_{\sL},$ and $\{\ph(\lsN)\:|\: \lN \in \Pi^n\}$ is a  basis of $\FL_{\sN}.$

Now we may describe the action of $\UN$ in the basis $\{\ph(\lsL)\:|\: \lL \in \Pi^l\}.$ 
For each $b\in \{1,\dots,l\}$ colour the nodes of the partition $\lambda^{(b)}$ by filling the node that lies at the intersection of $i$th row and $j$th column with $r=(s_b + j-i)\bmod n.$
 Write $\mL/\lL = \fra{r}$ to indicate that the multipartition $\mL$ is obtained from $\lL$ by adding a node with colour $r.$ 
In this case  $\mL/\lL $ is called a {\em removable $r$-node } of $\mL,$ or an {\em addable $r$-node} of $\lL.$ With every node of $\lL$ there is associated a pair $(d,b)$ of integers, where $d=s_b+j-i,$ and $b$ indicates the partition $\lambda^{(b)}$ to which the node belongs. 
A total order on the set of removable and addable nodes of $\lL$ is defined by:  
\begin{displaymath}
(d,b) < (d',b') \quad \Longleftrightarrow \quad ( (d < d') \; \text{or} \; (d=d' \: \text{and} \: b < b')).
\end{displaymath}
Now let $\mL/\lL$ be the $r$-node $(d,b).$ Put   
\begin{align*}
&N_r^<(\lL,\mL) =  \#\{ \text{addable $r$-nodes $(d',b')$ of $\lL$ such that $(d',b') < (d,b)$}\} \\
 &\qquad \qquad- \#\{ \text{removable $r$-nodes $(d',b')$ of $\lL$ such that $(d',b') < (d,b)$}\}, \\  
&N_r^>(\lL,\mL) =  \#\{ \text{addable $r$-nodes $(d',b')$ of $\lL$ such that $(d',b') > (d,b)$}\} \\
& \qquad \qquad- \#\{ \text{removable $r$-nodes $(d',b')$ of $\lL$ such that $(d',b') > (d,b)$}\}, \\
&N_r(\lL) = \#\{ \text{addable $r$-nodes  of $\lL$ }\} - \#\{ \text{removable $r$-nodes of $\lL$}\}.
\end{align*}

\begin{theor}
$\UN$ acts on $\FN_{\sL}$ by \\
\hspace{0.5cm} $\fN_r \ph(\lsL) = \sum_{\mL} q^{N_r^>(\lL,\mL)} \ph(\msL),$ sum over all $\mL$ such that  $\mL/\lL = \fra{r},$\\
\hspace{0.5cm} $\eN_r \ph(\msL) = \sum_{\lL} q^{-N_r^<(\lL,\mL)} \ph(\lsL),$ sum over all $\lL$ such that  $\mL/\lL = \fra{r},$\\
\hspace{0.5cm} $\tN_r \ph(\lsL) = q^{N_r(\lL)} \ph(\lsL).$  
\end{theor}
The $\UN$-module described by the above theorem was defined in \cite{JMMO}, where $\FN_{\sL}$ was called the (higher-level) $q$-deformed Fock space. 
The presentation of this action given above coincides with that of \cite{FLOTW} up to the transformation 
$$(\lambda^{(1)},\lambda^{(2)},\dots, \lambda^{(l)}) \mapsto (\lambda^{(l)},\dots, \lambda^{(2)},\lambda^{(1)}), \qquad  (s_1,s_2,\dots,s_l) \mapsto (s_l,\dots,s_2,s_1).$$ 


\section{Canonical bases of the $q$-deformed Fock spaces}

\subsection{The involution}
For a semi-infinite monomial (ordered or not)
\begin{displaymath}
u = u_{k_1}\wedge u_{k_2} \wedge \cdots \; \in \bigwedge^{s+\frac{\infty}{2}}\VNL
\end{displaymath}
set $\deg(u) := \sum_{i=1}^{\infty} k_i -(s-i+1).$ If $u \neq 0$ this is a non-negative integer, so $\deg :\bigwedge^{s+\frac{\infty}{2}}\VNL \rightarrow \Zint_{\geq 0}$ defines a grading of $\bigwedge^{s+\frac{\infty}{2}}\VNL.$ Let $r$ be a non-negative integer and let         
$$
C_r(u) = \sum_{1\leq i<j\leq r} \delta(b_i=b_j) - \delta(a_i=a_j),
$$
where, as usual, we put $k_i = a_i + n(b_i-1)-nlm_i$ ($a_i \in \{1,\dots,n\},$ $b_i\in \{1,\dots,l\},$ $m_i \in \Zint$).   

\begin{propos}
For  $r \geq \deg(u)$  the monomial
$$
\ov{u} = (-1)^{\frac{r(r-1)}{2}} q^{C_r(u)} u_{k_r}\wedge u_{k_{r-1}} \wedge \cdots \wedge u_{k_1} \wedge u_{k_{r+1}}\wedge u_{k_{r+2}} \wedge \cdots 
$$ 
is independent of $r.$ 
\end{propos}

Let $\{u^{(i)}\}$ be a collection of monomials, and let $c_i(q) \in \Qq.$ Define a semi-linear map $\barr : \bigwedge^{s+\frac{\infty}{2}}\VNL \rightarrow \bigwedge^{s+\frac{\infty}{2}}\VNL $ $:$ $v \mapsto \ov{v}$ by 
$$ 
\ov{ \sum_i c_i(q) u^{(i)} }  = \sum_i c_i(q^{-1}) \ov{u^{(i)}}. 
$$ 
\mbox{}From this definition it immediately follows that  $\barr$  is a degree-preserving involution of $\bigwedge^{s+\frac{\infty}{2}}\VNL.$  

\begin{theor} {\em (i)} For $v \in \bigwedge^{s+\frac{\infty}{2}}\VNL,$ 
$$
\ov{\fN_i(v)} = \fN_i(\ov{v}), \quad \ov{\fL_j(v)} = \fL_j(\ov{v}),\quad \ov{B_{-m}(v)}  = B_{-m}(\ov{v}).
$$
Here $i=0,1,\dots,n-1;$ $j=0,1,\dots,l-1;$ $m > 0.$ \\
{\em (ii)} For all $\sL =(s_1,\dots,s_l) \in \Zint^l$ $(\sum_{b=1}^l s_b= s)$ and all $\sN =(s_1,\dots,s_n) \in \Zint^n$ $(\sum_{a=1}^n s_a= s)$ the involution $\barr$ leaves invariant the  subspaces $\FN_{\sL}$ and $\FL_{\sN}.$ 
\end{theor}

Define matrices $\ANL_{\lambda,\mu}(s;q),$ $\ANL_{\lL,\mL}(\sL;q)$ and  $\ANL_{\lN,\mN}(\sN;q)$ by 
\begin{gather*}
\ov{\ph(\lambda,s)} = \sum_{\mu \in \Pi} \ANL_{\lambda,\mu}(s;q)\,\ph(\mu,s),\\ 
\ov{\ph(\lL,\sL)} = \sum_{\mL \in \Pi^l} \ANL_{\lL,\mL}(\sL;q)\,\ph(\mL,\sL) \quad \text{and} \quad 
\ov{\ph(\lN,\sN)} = \sum_{\mN \in \Pi^n} \ANL_{\lN,\mN}(\sN;q)\,\ph(\mN,\sN),
\end{gather*}
so that by (\ref{eq:ph1} -- \ref{eq:ph3}) :
\begin{alignat*}{2}
&\ANL_{\lL,\mL}(\sL;q) = \ANL_{\lambda,\mu}(s;q)& &\quad \text{for $(\lsL)=\iota_l(\lambda,s),$ $(\msL)=\iota_l(\mu,s),$} \\
&\ANL_{\lN,\mN}(\sN;q) = \ANL_{\lambda,\mu}(s;q)& &\quad \text{for $(\lsN)=\iota_n(\lambda,s),$ $(\msN)=\iota_n(\mu,s).$} 
\end{alignat*}

\begin{theor} 
{\em (i)} $ \ANL_{\lambda,\mu}(s;q) \in \Zint[q,q^{-1}].$\\
{\em (ii)} $\ANL_{\lambda,\lambda}(s;q) =1.$ \\  
{\em (iii)} $\ANL_{\lambda,\mu}(s;q) \neq 0$ only if $\lambda \geq \mu$ in the dominance ordering of partitions. \\ 
{\em (iv)} $\ALN_{\lN,\mN}(\sN;q) = \ANL_{\lN,\mN}(\sN;-q^{-1}),$ and $\ALN_{\lL,\mL}(\sL;q) = \ANL_{\lL,\mL}(\sL;-q^{-1}).$ 
\end{theor}

For a multipartition $\lL = (\lambda^{(1)},\dots,\lambda^{(l)})$ set $|\lL| = \sum_{b=1}^l |\lambda^{(b)}|.$ 
\mbox{} From the  fact that $\barr$ is degree-preserving, and Theorem 2.2 (ii) it follows  that $\ANL_{\lL,\mL}(\sL;q) \neq 0$ only if $|\lL| = |\mL|,$ and that  $\ANL_{\lN,\mN}(\sN;q) \neq 0$ only if $|\lN| = |\mN|.$ Therefore the vectors $\ph(\bo_l,\sL)$ and $\ph(\bo_n,\sN)$ are both $\barr$-- invariant.  

By Theorem 1.1 the subspace $M_{\sL} = \UN\cdot\ph(\bo_l,\sL)\subset \FN_{\sL}$ is isomorphic to the irreducible $\UN$-module $V(\Lambda)$ with highest weight $\Lambda = \Lambda_{s_1\bmod n} + \cdots + \Lambda_{s_l\bmod n}.$ 
Now Theorem 2.2 (i) and $ \ov{\ph(\bo_l,\sL)} = \ph(\bo_l,\sL)$ imply that  $M_{\sL}$ is $\barr$ -- invariant, and that the restriction of $\barr$ on $M_{\sL}$ coincides, under the above isomorphism, with the standard involution in terms of which global crystal bases of  $V(\Lambda)$  are defined \cite{K}.    

\subsection{Canonical bases}

For $s\in \Zint$ let $L^{\pm}(s) = {\mathbb Q}[q^{\pm 1}] \{ \ph(\lambda,s)\:|\: \lambda \in \Pi\}.$ 
\begin{theor}
There are unique bases $\{G^+(\lambda,s)\: | \: \lambda \in \Pi\},$ $\{G^-(\lambda,s)\: | \: \lambda \in \Pi\}$ of $\bigwedge^{s+\frac{\infty}{2}}\VNL$ such that:\\  
{\em (i)} $ \ov{G^{\pm}(\lambda,s)} = G^{\pm}(\lambda,s),$ \\ 
{\em (ii)} $ G^{\pm}(\lambda,s) \equiv \ph(\lambda,s) \bmod q^{\pm 1} L^{\pm}(s).$
\end{theor}
Define 
\begin{alignat}{4}
&G^{\pm}(\lsL)& &:= &  &\;G^{\pm}(\lambda,s)& &\qquad\text{ for $(\lsL) = \iota_l(\lambda,s),$}\\
&G^{\pm}(\lsN)& &:= &  &\;G^{\pm}(\lambda,s)& &\qquad\text{ for $(\lsN) = \iota_n(\lambda,s).$} 
\end{alignat}
Then $\{G^+(\lsL)\: | \: \lL \in \Pi^l\},$ $\{G^-(\lsL)\: | \: \lL \in \Pi^l\}$ are bases of $\FN_{\sL}.$ We  call the $\{G^+(\lsL)\}$ the {\em canonical basis} of the $q$-Fock space  $\FN_{\sL}.$ 

Now define matrices $\DNLpm_{\lambda,\mu}(s;q),$ $\DNLpm_{\lL,\mL}(\sL;q)$ and $\DNLpm_{\lN,\mN}(\sN;q)$ by  
\begin{gather*}
G^{\pm}(\lambda,s) = \sum_{\mu \in \Pi}\DNLpm_{\lambda,\mu}(s;q) \, \ph(\mu,s),  \mbox{} \\
G^{\pm}(\lsL) = \sum_{\mL \in \Pi^l} \DNLpm_{\lL,\mL}(\sL;q)\, \ph(\msL) \;\; \text{and}\;\;   
G^{\pm}(\lsN) = \sum_{\mN \in \Pi^n} \DNLpm_{\lN,\mN}(\sN;q)\, \ph(\msN). 
\end{gather*}
By (\ref{eq:ph1} -- \ref{eq:ph3}) we have 
\begin{alignat*}{2}
&\DNLpm_{\lL,\mL}(\sL;q) = \DNLpm_{\lambda,\mu}(s;q)& &\quad \text{for $(\lsL)=\iota_l(\lambda,s),$ $(\msL)=\iota_l(\mu,s),$} \\
&\DNLpm_{\lN,\mN}(\sN;q) = \DNLpm_{\lambda,\mu}(s;q)& &\quad \text{for $(\lsN)=\iota_n(\lambda,s),$ $(\msN)=\iota_n(\mu,s).$} 
\end{alignat*}
Theorem 2.3 implies that $\DNLpm_{\lambda,\mu}(s;q) \in \Zint[q^{\pm 1}],$  $\DNLpm_{\lambda,\mu}(s;q) \neq 0$ only if $\lambda \geq \mu,$ and   
\begin{alignat*}{3}
&\DLNpm_{\lL,\mL}(\sL;q)&  &=  & &\;\DNLmp_{\lL,\mL}(\sL;-q^{-1}), \\
&\DLNpm_{\lN,\mN}(\sN;q)&  &= & &\;\DNLmp_{\lN,\mN}(\sN;-q^{-1}). 
\end{alignat*}

Let $n>s_1\geq \cdots \geq s_l \geq 0.$ Following \cite{FLOTW} we call $\lL \in \Pi^l$ $\sL$-{\em cylindrical} if  
\begin{alignat*}{3}
&\lambda^{(b+1)}_i &  &\geq  \lambda^{(b)}_{i+s_b-s_{b+1}} & &\quad 1\leq b < l,\; i=1,2,\dots,\\
&\lambda^{(1)}_i   &  &\geq \lambda^{(l)}_{i+n+s_l-s_{1}}  & &\quad  i=1,2,\dots.
\end{alignat*}
Let $\Pi^l_{\sL}$ be the set of all $\sL$-{cylindrical} multipartitions $\lL$ such that for all $k>0,$ among colours appearing at the right ends of the length $k$ rows of $\lL$ at least one element of $\{0,1,\dots,n-1\}$ does not occur. It follows from \cite[Proposition 2.11]{FLOTW} that $\{G^+(\lsL)\: | \: \lL \in \Pi^l_{\sL}\} $ is the global lower crystal basis of the irreducible submodule $M_{\sL}$ of $\FN_{\sL}.$

\subsection{Canonical bases and Kazhdan-Lusztig polynomials}
In this section we give expressions for the transition matrices $\DNLpm_{\lambda,\mu}(s;q)$ in terms of (parabolic) affine Kazhdan-Lusztig polynomials. Our notations concerning the affine Weyl group and the affine Hecke algebra mostly follow the work \cite{LT2}. 
\subsubsection{Affine Weyl group}
Let $P= \oplus_{i=1}^r \Zint \epsilon_i$ be the weight lattice of ${\mathfrak {gl}}_r.$ The Weyl group $W = {\mathfrak S}_r$ with the system of generators $s_i$ $(i=1,\dots,r-1)$ acts on $P$ from the left with $s_i$ exchanging $\epsilon_i $ and $\epsilon_{i+1}.$ 
The affine Weyl group $\widehat{W}$ is the semi-direct product $W\ltimes P.$ Let $\pi = \epsilon_1 s_1\cdots s_{r-1},$ $s_0 = \pi s_{r-1} \pi^{-1},$ then $\pi^{\pm 1},s_i$ $(i=0,1,\dots,r-1)$ is a system of generators of $\widehat{W}.$   

The group $\widehat{W}$ acts on $\Zint^r$ from the right by
\begin{align*}
&(h_1,\dots,h_r)\cdot s_i = (h_1,\dots,h_{i+1},h_i,\dots,h_r) \qquad (i=1.\dots,r-1),\\
&(h_1,\dots,h_r)\cdot \epsilon_i = (h_1,\dots,h_i+n,\dots h_r)\qquad (i=1.\dots,r).
\end{align*}
The set $\A_r^n = \{ (a_1,\dots,a_r) \:|\: 1\leq a_1\leq \cdots \leq a_r \leq n\}$ is a fundamental domain of this action. For $\ba \in \A_r^n$ denote by $W_{\ba}$ the stabilizer of $\ba,$ and let $\aW$ be the set of minimal length representatives in $W_{\ba}\!\setminus\!\widehat{W}.$ 

\subsubsection{Affine Hecke algebra}

The Hecke algebra $\widehat{H}$ is the algebra over $\Zint[q,q^{-1}]$ with basis $\{ T_x\:|\: x \in \widehat{W}\}$ and relations 
\begin{align*}
&T_x T_y = T_{xy} \quad \text{if} \quad l(x)+l(y) = l(xy),\\
&(T_{s_i} - q^{-1}) (T_{s_i} + q) = 0 \quad (i=0,1,\dots,r-1).
\end{align*}
There is a canonical involution $x \mapsto \ov{x}$ of $\widehat{H}$ defined as the unique ring homomorphism such that $\ov{q} = q^{-1}$ and $\ov{T_x} = (T_{x^{-1}})^{-1}.$ The Hecke algebra has unique bases $\{ C_x \: | \: x \in \widehat{W}\}$ and $\{ C_x' \: | \: x \in \widehat{W}\}$ characterized by   
\begin{gather*}
 \ov{C_x} =  C_x,\qquad  \ov{C_x'} =  C_x', \\ 
C_x \equiv T_x \bmod q^{-1}\Zint[q^{-1}]\{ T_y \}, \qquad C_x' \equiv T_x \bmod q\Zint[q]\{ T_y \}. 
\end{gather*}
Put
$$ 
C_x = \sum_{y} \P_{y,x}^- \,T_y,\qquad C_x' = \sum_{y} \P_{y,x}^+ \,T_y.
$$
Then $\P_{y,x}^+ = q^{l(x)-l(y)}P_{y,x},$ $\P_{y,x}^- = (-q)^{l(y)-l(x)}\ov{P_{y,x}},$ where $P_{y,x} \in {\mathbb N}[q^{-2}]$ are the Kazhdan-Lusztig polynomials.  
\vspace{3mm}

\noindent To express the matrices  $\DNLpm_{\lambda,\mu}(s;q)$ in terms of $\P^{\pm}_{y,x}$ we need to prepare some notations. With a semi-infinite integral sequence $\bk =(k_1,k_2,\dots)$ such that
 $k_1>k_2>\cdots,$ $k_i = s-i+1$ $(i\gg 1),$ we associate $\ba(\bk) \in \A_r^n,$ $x(\bk) \in \akW,$ $\bb(\bk) \in \B_r^l = \{(b_1,\dots,b_r) \in \Zint^r\: |\: l\geq b_1\geq \cdots \geq b_r \geq 1\},$ where $r=\deg(u_{k_1}\wedge u_{k_2}\wedge \cdots\:) $ 
$=$ $\sum_{i\geq 1} k_i-s+i-1,$ as follows. Put, as usual, $k_i = a_i +n(b_i -1) -nl m_i$ ($a_i \in \{1,\dots,n\},$ $b_i\in \{1,\dots,l\},$ $m_i \in \Zint$). 

Now define $\ba(\bk)$ to be $(a_1,\dots,a_r)$ rearranged in the non-decreasing order of magnitude, and define $\bb(\bk)$ to be $(b_1,\dots,b_r)$ rearranged in the non-increasing order of magnitude. 

For each $b=1,2,\dots,l$ let $\bk^{(b)}$ be the integer vector obtained by arranging elements of the set $\{ a_i -nm_i\:|\: b_i =b,\, 1\leq i\leq r\}$ in the strictly decreasing order. Put $\bh = (\bk^{(l)},\bk^{(l-1)},\dots,\bk^{(1)}).$ Then $\bh \in \Zint^r,$ and $\bh \in \ba(\bk)\cdot \widehat{W}.$  
Now define $x(\bk)$ to be the preimage of $\bh$ under the bijection $$\akW \rightarrow \ba(\bk)\cdot \widehat{W}\quad :\quad x \mapsto \ba(\bk) \cdot x.$$  

Let $\DNLpm_{\bk, \bll} = \DNLpm_{\lambda,\mu}(s;q)$ where $\bk = (k_1,k_2,\dots\;),$ $\bll = (l_1,l_2,\dots\;),$ and $k_i = \lambda_i + s-i+1,$ $l_i = \mu_i + s-i+1.$ Observe that Theorem 2.2 (ii) 
implies that $\DNLpm_{\bk, \bll} \neq 0$ only if $\ba(\bk) = \ba(\bll),$ $\bb(\bk) = \bb(\bll).$ Set $\ba = \ba(\bk), \bb= \bb(\bk),$ and let $W_{\ba},W_{\bb} \subset W$ be the stabilizers of $\ba,\bb$ respectively. 
Denote by $\omega_{\ba},\omega_{\bb}$ the longest elements in $W_{\ba},W_{\bb}.$   

\begin{theor}
\begin{alignat*}{2}
&\DNLm_{\bk, \bll}& &= \sum_{\sigma \in W_{\ba}} q^{-l(\sigma)} \P^-_{\sigma x(\bll),\,x(\bk)}, \\ 
&\DNLp_{\bk, \bll}& &= \sum_{\sigma \in W_{\bb}} (-q)^{l(\sigma)} \P^+_{ \omega_{\ba}x(\bll)\omega_{\bb}\sigma ,\,\omega_{\ba}x(\bk)\omega_{\bb}}. \\ 
\end{alignat*}
\end{theor}
\noindent A proof of this theorem, as well as proofs of other results of this note will be given in a future publication. 
\section{Appendix}

In this appendix we give examples of the matrices $\DNLp_{\lL,\mL}(\sL;q)$ for $n=l=2,$ $\sL = (0,0).$ The rows of the matrices are labelled by  $ \mu | \mL,$ where $(\msL) = \iota_l(\mu,0).$ 
The matrices are given up to $|\lL | = 4,$ and are split into blocks according to the value of $|\lambda|.$ For example, the entry corresponding to $\lL= ((4),\emptyset),$ $\mL= ((2,1),(1))$ equals $2q^2.$

\begin{small}

\noindent $|\lL  | = 1$
%
\begin{displaymath}
\begin{array}{l|l| c c}
(3) &  (\emptyset, (1)) &    1 & \cdot \cr 
(1^3) &  ((1), \emptyset) &  q & 1    
\end{array}
\end{displaymath}
$|\lL  | = 2$
%
\begin{displaymath}
\begin{array}{l|l| c c c c }
(4) &  (\emptyset, (2)) &  1 & \cdot & \cdot & \cdot \cr 
(3, 1) &  (\emptyset, (1^2)) & q & 1 & \cdot & \cdot \cr 
(2, 1^2) &  ((2), \emptyset) & q & \cdot & 1 & \cdot \cr 
(1^4) &  ((1^2), \emptyset) & {q^2} & q & q & 1 
\end{array}\qquad 
\begin{array}{l|l| c } 
 (3, 2, 1) &  ((1), (1)) & 1  
\end{array}
\end{displaymath} 
%
%
%
%
%
%
%
%
%
%
$|\lL  | = 3$
\begin{displaymath}
\begin{array}{l|l| c c}
(4, 1) &  (\emptyset, (2, 1)) & 1 & \cdot \cr 
(2, 1^3) &  ((2, 1), \emptyset) & q & 1    
\end{array}\qquad
\begin{array}{l|l| c c c c c c c c } 
(7) &  (\emptyset, (3)) &   1 & \cdot & \cdot & \cdot & \cdot & \cdot & \cdot & \cdot \cr 
(5, 1^2) &  ((3), \emptyset) & q & 1 & \cdot & \cdot & \cdot & \cdot & \cdot & \cdot \cr 
(4, 2, 1) &  ((1), (2)) & q & \cdot & 1 & \cdot & \cdot & \cdot & \cdot & \cdot \cr 
(3^2 , 1) &  ((2), (1)) & {q^2} & q & q & 1 & \cdot & \cdot & \cdot & \cdot \cr 
(3, 2^2) &  ((1), (1^2)) & q & \cdot & {q^2} & q & 1 & \cdot & \cdot & \cdot \cr 
(3, 2, 1^2) &  ((1^2), (1)) & {q^2} & q & {q^3} & {q^2} & q & 1 & \cdot & \cdot \cr 
(3, 1^4) &  (\emptyset, (1^3)) & {q^2} & \cdot & \cdot & \cdot & q & \cdot & 1 & \cdot \cr 
(1^7) &  ((1^3), \emptyset) & {q^3} & {q^2} & \cdot & q & {q^2} & q & q & 1  
\end{array}
\end{displaymath}
$|\lL | = 4$
\begin{displaymath}
\begin{array}{l | l |   c c c c c c c c c c c c c c c c}
(8) &  (\emptyset, (4))  &    
1 & \cdot  & \cdot  & \cdot  & \cdot  & \cdot  & \cdot  & \cdot  & \cdot  & \cdot  & \cdot  & \cdot  & \cdot  & \cdot  & \cdot  & \cdot  \cr 
(7, 1)& (\emptyset, (3, 1))  &  
q & 1 & \cdot  & \cdot  & \cdot  & \cdot  & \cdot  & \cdot  & \cdot  & \cdot  & \cdot  & \cdot  & \cdot  & \cdot  & \cdot  & \cdot  \cr 
(6, 1^2) &  ((4), \emptyset)  &  
q & \cdot  & 1 & \cdot  & \cdot  & \cdot  & \cdot  & \cdot  & \cdot  & \cdot  & \cdot  & \cdot  & \cdot  & \cdot  & \cdot  & \cdot  \cr 
(5, 1^3) &  ((3, 1), \emptyset)  &  
{q^ 2} & q & q & 1 & \cdot  & \cdot  & \cdot  & \cdot  & \cdot  & \cdot  & \cdot  & \cdot  & \cdot  & \cdot  & \cdot  & \cdot  \cr 
(4^2) &  (\emptyset, (2^2))  &   
\cdot  & q & \cdot  & \cdot  & 1 & \cdot  & \cdot  & \cdot  & \cdot  & \cdot  & \cdot  & \cdot  & \cdot  & \cdot  & \cdot  & \cdot  \cr 
(4, 3, 1) &  ((2), (2))  &      
{q^2} & \cdot  & q & \cdot  & \cdot  & 1 & \cdot  & \cdot  & \cdot  & \cdot  & \cdot  & \cdot  & \cdot  & \cdot  & \cdot  & \cdot  \cr 
(4, 2^2) &  ((1), (2, 1))  &    
q & {q^2} & \cdot  & \cdot  & q & q & 1 & \cdot  & \cdot  & \cdot  & \cdot  & \cdot  & \cdot  & \cdot  & \cdot  & \cdot  \cr 
(4, 2, 1^2) &  ((1^2), (2))  &  
{q^2} & \cdot  & q & \cdot  & \cdot  & {q^ 2} & q & 1 & \cdot  & \cdot  & \cdot  & \cdot  & \cdot  & \cdot  & \cdot  & \cdot  \cr 
(4, 1^4) &  (\emptyset, (2, 1^2))  & 
{q^2} & {q^ 3} & \cdot  & \cdot  & {q^2} & \cdot  & q & \cdot  & 1 & \cdot  & \cdot  & \cdot  & \cdot  & \cdot  & \cdot  & \cdot  \cr 
(3^2, 2) &  ((2), (1^2))  &       
{q^ 2} & \cdot  & q & \cdot  & \cdot  & {q^2} & q & \cdot  & \cdot  & 1 & \cdot  & \cdot  & \cdot  & \cdot  & \cdot  & \cdot  \cr 
(3^2, 1^2) &  ((2, 1), (1))  &     
{q^3} & {q^2} & 2\,{q^2} & q & q & {q^3} & {q^2} & q & \cdot  & q & 1 & \cdot  & \cdot  & \cdot  & \cdot  & \cdot  \cr 
(3, 2^2, 1) &  ((1^2), (1^2))  &     
{q^ 2} & \cdot  & {q^3} & \cdot  & {q^2} & {q^4} & q + {q^3} & {q^2} & \cdot  & {q^2} & q & 1 & \cdot  & \cdot  & \cdot  & \cdot  \cr 
(3, 1^5) &  (\emptyset, (1^4))  &     
{q^3} & \cdot  & \cdot  & \cdot  & \cdot  & \cdot  & {q^2} & \cdot  & q & \cdot  & \cdot  & q & 1 & \cdot  & \cdot  & \cdot  \cr 
(2^4) &  ((2^2), \emptyset)  &          
\cdot  & {q^3} & \cdot  & {q^2} & {q^2} & \cdot  & \cdot  & \cdot  & \cdot  & \cdot  & q & \cdot  & \cdot  & 1 & \cdot  & \cdot  \cr 
(2, 1^6) &  ((2, 1^2), \emptyset)  &     
{q^3} & {q^4} & {q^2} & {q^3} & {q^3} & \cdot  & {q^2} & q & q & q & {q^2} & \cdot  & \cdot  & q & 1 & \cdot  \cr 
(1^8) &  ((1^4), \emptyset)  &         
{q^4} & \cdot  & {q^3} & \cdot  & \cdot  & \cdot  & {q^3} & {q^2} & {q^2} & {q^2} & q & {q^2} & q & \cdot  & q & 1  
\end{array}
\end{displaymath}

%
\begin{displaymath}
\begin{array}{l|l| c c c c  }
(7, 2, 1) &  ((1), (3)) &  1 & \cdot &  \cdot & \cdot \cr 
(5, 4, 1) &  ((3), (1)) &  q & 1 &  \cdot & \cdot \cr 
(3, 2^3, 1) &  ((1), (1^3))  & q & \cdot &  1 & \cdot \cr 
(3, 2, 1^5) &  ((1^3), (1)) &  {q^2} & q &  q & 1  
\end{array}
\end{displaymath}
\end{small}
\newcommand{\BOOK}[6]{\bibitem[{#6}]{#1}{\sc #2}, {\it #3} (#4)#5.}
\newcommand{\JPAPER}[8]{\bibitem[{#8}]{#1}{\sc #2}, {\it #3},
{#4} {\bf #5} (#6) #7.}
\newcommand{\PRE}[6]{\bibitem[{#6}]{#1}{\sc #2}, {\it #3},
{#4}  (#5).}
\newcommand{\JPAPERS}[9]{\bibitem[{\bf #9}]{#1}{\sc #2}, `#3', {\it #4} #5 #6,
#7 #8.}
%

%
%
\begin{it}
Research Institute for Mathematical Sciences
\hfill
E-mail\::\hspace{1mm}
{\em duglov@kurims.kyoto-u.ac.jp}

\noindent
Kyoto University

\noindent
Kyoto 606-8502, JAPAN
\end{it}
\vspace{-5mm}
\end{document}